\documentclass[12pt]{article}

\title{The area of a self-similar fragmentation} 
\author{Jean Bertoin\thanks{Laboratoire de Probabilit\'es et Mod\`eles Al\'eatoires, 
UPMC, 4 Place Jussieu, 75252 Paris Cedex 05; France. 
Email: jean.bertoin@upmc.fr} 
 }

\date{}
\usepackage{amsfonts,amsmath,amssymb,bbm}
\usepackage{setspace} 
\def\proof{\noindent{\bf Proof:}\hskip10pt}        
\def\QED{\hfill $\Box$}

\font\tenmath=msbm10 scaled 1200
\font\sevenmath=msbm7 scaled 1200
\font\Fivemath=msbm5 scaled 1200
\newfam\mathfam \textfont\mathfam=\tenmath
\scriptfont\mathfam=\sevenmath \scriptscriptfont\mathfam=\Fivemath

\def \\ { \cr }
\def\R{\mathbb{R}}
\def \1{1 \mkern -6mu 1} 
\def\N{\mathbb{N}}
\def\E{\mathbb{E}}
\def\P{\mathbb{P}}

\def\R{\mathbb{R}}

\def \n{^{(n)}}

\def \e{{\rm e}}
\def \x{{\bf x}}
\def \X{{\bf X}}

\def \d{{\rm d}}
\def \p{{\mathcal P}}

\def \a{{\mathcal A}}
\def \m{\p_{\rm m}}

\newtheorem{theorem}{Theorem}

\newtheorem{proposition}{Proposition}
\newtheorem{lemma}{Lemma}
\newtheorem{corollary}{Corollary}
\begin{document}

\maketitle

\begin{abstract} 
We consider the area $A=\int_0^{\infty}\left(\sum_{i=1}^{\infty} X_i(t)\right) \d t$ 
of a self-similar fragmentation process $\X=(\X(t), t\geq 0)$ with negative index.
We  characterize the law of $A$ by an integro-differential equation. The latter may be viewed as
the infinitesimal version of a recursive distribution equation
that arises naturally in this setting.  In the case of binary splitting, this yields a recursive formula 
for the entire moments of $A$ which
generalizes  known results 
for the area of the Brownian excursion. 

\end{abstract}

{\bf Key words:} Self-similar fragmentation, area, recursive distributional equation.
  
\begin{section}{Introduction}

The distribution of the area  $A_{\rm Exc}=\int_0^1e_s \d s$  of a standard Brownian excursion
$(e_s, 0\leq s \leq 1)$ appears in a variety of settings, including random graphs
\cite{spencer, tak1}, random trees and branching processes \cite{tak2, tak4}, order statistics \cite{tak3}, hashing with linear probing \cite{FPV}, ..., not to mention of course the study of Brownian motion for its own interest \cite{louchard1, louchard2} . 
The entire moments $\E( A_{\rm Exc}^k)$ have a special importance, as they are related, for instance, to asymptotics as $n\to\infty$ for the number of connected graphs with $n$ labelled vertices and $n+k-1$ edges, see \cite{spencer} and the survey  \cite{janson}. We refer to Perman and Wellner \cite{PW}, Janson \cite{janson} and references therein for a detailed presentation  and reviews of known results on this topic.

The starting point of this work lies in the observation that one can express the area in the form
$$A_{\rm Exc}=\int_0^{\infty} |\theta(t)| \d t$$
where $ |\theta(t)|$ denotes the Lebesgue measure of the random open set $\theta(t)=\{s\in [0,1]: e_s>t\}$.
The point is that the process $\theta=( \theta(t), t\geq 0)$ is a self-similar interval-fragmentation
 in the sense of \cite{Be1, RFCP}. One can  derive 
an integro-differential equation for the distribution of $A_{\rm Exc}$ from the branching and self-similarity properties of $\theta$. In particular this yields
recursive formulas for the entire moments of $A_{\rm Exc}$ that have been obtained in the literature by analytic techniques based on the Feynman-Kac formula or the analysis of continued fractions and of singularities of  the generating functions of discrete approximations.

The same approach applies more generally to self-similar fragmentation processes, a class
of Feller processes  with values in the space of
mass-partitions 
$$\m=\{\x=(x_1, x_2, \ldots): x_1\geq x_2\geq \ldots \geq 0\hbox{ and }
\sum_{i=1}^{\infty}x_i\leq 1\}\,.$$
Specifically, a self-similar fragmentation  $\X=(\X(t), t\geq 0)$
 fulfills the following two fundamental properties. For $x\in[0,1]$, let us denote by $\P_x$ the law of the version of $\X$ started from a single fragment of mass $x$, i.e.  $\X(0)=(x, 0, \ldots)$. First, the self-similarity means 
that there exists an index $\alpha\in\R$ such that for every $x\in(0,1]$ the distribution
of the rescaled process $(x\X(x^{\alpha}t), t\geq 0)$ under $\P_1$ is $\P_x$. Second, the process $\X$ satisfies the branching property, in the sense that for every mass-partition $\x=(x_1, x_2, \ldots)$,
if $\X^{(1)}, \X^{(2)}, \ldots$ are independent fragmentations with respective laws
$\P_{x_1}, \P_{x_2}, \ldots$, then the process resulting from the decreasing rearrangement
of all the fragments of $\X^{(1)}(t), \X^{(2)}(t), \ldots$ is a version of $\X(t)$ started from $\X(0)=\x$.

We assume that the index of self-similarity $\alpha$ is negative, which implies that small fragments 
split faster than the large ones. A well-known consequence is that
the process of the total mass $t\mapsto \sum_{i=1}^{\infty}X_i(t)$ decreases  
and reaches $0$ in finite time a.s.; in other words the entire mass is eventually ground down to dust. This has been observed first  by Filippov \cite{Fi}, see also Proposition 2(i) in \cite{Be2}. We may thus define the area
$$A=\int_0^{\infty}\left(\sum_{i=1}^{\infty} X_i(t)\right) \d t$$
which is the main object of interest in this work.
The denomination {\it area} is better understood if we remember that a fragmentation admits an interval representation; cf. Section 3.2 in \cite{Be1}. There exists a nested right-continuous family $(\theta(t), t\geq 0)$ of open subsets of
the unit interval such that for every $t\geq 0$, the sequence of the lengths of the interval components of $\theta(t)$ listed in the decreasing order is precisely 
$\X(t)= (X_1(t), \ldots)$. 
If we define a lower semi-continuous path $F:[0,1]\to\R_+$ by
$$F(u)=\sup\{t\geq 0: u\in \theta(t)\}\,,\qquad u\in [0,1]\,,$$
then $\theta(t)=\{u: F(u)>t\}$, and since $\sum_{i=1}^{\infty} X_i(t)=|\theta(t)|$, we can express $A$ in the form 
$A=\int_0^1F(u)\d u$. 
Of course $F=e$ is the Brownian excursion in the situation discussed at the beginning of this introduction.

The area $A$ has another natural interpretation in terms of continuous random trees.
Indeed, Haas and Miermont \cite{HM} obtained a representation of self-similar fragmentations in terms of some rooted continuous tree ${\bf T}$ which enjoys a self-similarity and branching  properties.
More precisely,  for every $t\geq 0$, $\X(t)= (X_1(t), \ldots)$
can be viewed as the ranked sequence of the masses of connected components of
${\bf T}(t)$,  the subset
of  the points in ${\bf T}$ at distance at least $t$ from the root.  In this setting,
$$A=\int_0^{\infty}|{\bf T}(t)| \d t \,,$$
where $|{\bf T}(t)|$ denotes the mass of ${\bf T}(t)$. Hence  $A$ represents the average 
height in ${\bf T}$, i.e. the average distance
of points to the root, where averaging is taken with respect to the mass measure of ${\bf T}$.
We also refer to \cite{tak4} for results on this quantity in the framework of certain discrete random trees.

In the next section,  we will present our main result which determines  the law of $A$ as the unique solution to an intro-differential equation expressed in terms of the characteristics of $\X$. 
In the case of binary dislocations, this enables us to derive explicit recursive formulas for the moments
of $A$. We recover in particular identities due to Tak\'acs \cite{tak1}
for the moments of the area of the Brownian excursion. 
Section 3 is devoted to the proof of  this integro-differential equation. We shall start by establishing 
a priori bounds for the moments of $A$. Then we proceed with the simpler case 
when the dislocation measure is finite, and derive the equation from a recursive distribution equation which is naturally induced by  the dynamics of fragmentation. The general case when the dislocation measure is infinite 
is then deduced by approximation. This relies on a weak limit theorem for the area of fragmentation processes. The proof of the latter is somewhat technical and will be postponed to the final subsection.
\end{section}

\begin{section}{Main results}

We denote by $\nu$ the dislocation measure of $\X$, so $\nu$ is a sigma-finite measure
on $\m$ with no atom at the trivial mass-partition $(1, 0, \ldots)$ and fulfills the integral condition
\begin{equation}\label{E9}
\int_{\m}(1-x_1)\nu(\d \x)<\infty\,.
\end{equation}
We implicitely exclude the degenerate case when $\nu\equiv 0$ and further assume absence of erosion. Roughly speaking, this means that $\X$ is a purely discontinuous process that only evolves by sudden dislocations whose rates are determined by $\nu$ and the index of self-similarity $\alpha$. We refer to Chapter 3 of \cite{RFCP} or \cite{Be1} for background.

For an arbitrary mass partition $\x$, we denote by 
$$\eta_{\x}(\d a) = \P_{\x}(A\in \d a)\,,\qquad a\geq 0\,, $$
the law of the area under the probability measure $\P_{\x}$
for  which the fragmentation $\X$ starts from $\X(0)=\x$.
For the sake of simplicity, we will work from now on under the law $\P=\P_1$, i.e. when the fragmentation starts from a single fragment of mass $1$, and write then $\eta=\eta_1$.
 This induces no loss of generality since, combining self-similarity and the branching property, we get that for every mass-partition $\x=(x_1, x_2, \ldots)$,  there is 
 the identity
 \begin{equation}\label{E1}
\eta_{\x}(\d a) = \P\left(\sum_{i=1}^{\infty} x_i^{1-\alpha} A_i\in \d a\right)
\end{equation}
where $(A_i)_{i\in\N}$ is a family of i.i.d. copies of $A$.
Note that when $\x$ has only finitely many non-zero terms, say $x_1,\ldots, x_n$, then
$\eta_{\x}$ can be expressed as a convolution product
$\eta_{\x}=\eta_{x_1}*\cdots * \eta_{x_n}$ where $\eta_y$ stands for the image of $\eta$
by the dilation $a\mapsto y^{1-\alpha}a$.
Finally, we let $\langle \mu, f\rangle=\int f \d \mu$ denote the integral of some function $f$ with respect to a measure $\mu$ when the integral makes sense. We are now able to state our main result

\begin{theorem}\label{T1} Let $\X$  be a self-similar fragmentation with index $\alpha<0$,  
dislocation measure $\nu$ and without erosion.
Then for every ${\mathcal C}^1$ function $f: \R_+\to \R$ such that $f'(y)=O(y^p)$ as $y\to\infty$  for some $p>0$, 
the law $\eta$  of $A$ solves
\begin{equation}\label{E0}
\langle \eta, f'\rangle = \int_{\m} \nu(\d \x) \left( \langle \eta, f\rangle - \langle \eta_{\x}, f\rangle \right)\,,
\end{equation}
where the quantities above are finite.
Further \eqref{E0} characterizes $\eta$.
\end{theorem}

If we introduce the concave increasing function $\Phi: \R_+\to \R_+$ by
\begin{equation}\label{E8}
\Phi(q)=\int_{\m}\left(1-\sum_{j=1}^{\infty}x_j^{1+q}\right)\nu(\d \x)\,,
\end{equation}
then we immediately see from Theorem \ref{T1}  that the first moment of $A$ is given simply by
$$\E(A)=\langle \eta, {\rm Id}\rangle=1/\Phi(-\alpha)\,.$$
This identity  can also be established directly; see the forthcoming Lemma \ref{L1}.

More generally, we shall now derive from Theorem \ref{T1} a recursive formula for the entire moments of $A$.
For the sake of simplicity, we shall focus on the special case of binary dislocations, although
more general situations could be dealt with at the price of heavier notation. This means that we assume  that the dislocation measure $\nu$ has support in the subset of binary mass-partitions $\{\x=(x,1-x,0, \ldots), x\in[1/2,1)\}$. By a slight abuse, 
we shall then identify the dislocation measure $\nu$ with its image by the map $\x\to x_1$,
i.e. we view $\nu$ as a measure on $[1/2,1)$. Specializing Theorem \ref{T1} to $f(x)=x^k$ and
applying the binomial formula, we immediately obtain:

\begin{corollary} \label{C1} Let $\X$  be a self-similar fragmentation with index $\alpha<0$,  
binary dislocation measure $\nu$ and without erosion. 
For every integer $k\geq 0$, let $M_k=\E(A^k)$ denote the $k$-th moment of the area.
Then there is the identity
$$a_k M_k = kM_{k-1} + \sum_{j=1}^{k-1} a_{j,k} M_jM_{k-j}\,, \qquad k\geq 1\,,$$
with
$$a_k=\int_{[1/2,1)}(1-x^{k(1-\alpha)}-(1-x)^{k(1-\alpha)})\nu(\d x)=\Phi(k(1-\alpha)-1)\,,$$
and
$$a_{j,k}=\left(\begin{matrix}
k\\j \end{matrix} \right)\int_{[1/2,1)}x^{j(1-\alpha)}(1-x)^{(k-j)(1-\alpha)}\nu(\d x)\,.$$
 \end{corollary} 
We stress that \eqref{E9} ensures the finiteness of $a_k$ and $a_{j,k}$. 

We now discuss some examples, starting with 
 the case of the Brownian fragmentation $A=A_{\rm Exc}$ which has motivated this work. 
 The Brownian fragmentation has self-similarity index  $\alpha=-1/2$, no erosion, and its dislocation measure is binary  and specified by
$$\nu(\d x) = \frac{2}{\sqrt{2\pi x^3 (1-x)^3}}\d x\,,\qquad 1/2\leq x < 1\,;$$
see \cite{Be1} on its pages 339-340. One gets by symmetry
$$
a_{k}\,=\, \int_0^1 \frac{1- x^{3k/2}-(1-x)^{3k/2}}{\sqrt{2\pi x^3 (1-x)^3}}\d x\,=\, 2^{3/2} \frac{\Gamma((3k-1)/2)}{\Gamma(3k/2-1)}
$$
and 
\begin{eqnarray*}
a_{j,k}&=& \left(\begin{matrix}
k\\j \end{matrix} \right)\int_0^1 \frac{x^{3j/2}(1-x)^{3(k-j)/2}}{\sqrt{2\pi x^3 (1-x)^3}}\d x\\
&=& \frac{k!\, \Gamma((3j-1)/2)\,\Gamma((3(k-j)-1)/2)}{\sqrt{2\pi}\, j! \, (k-j)! \,\Gamma(3k/2-1)}
\end{eqnarray*}
Following Tak\'acs \cite{tak1}, if we set
$$M_k=\frac{4 \sqrt{\pi} 2^{-k/2} k!}{\Gamma((3k-1)/2)} K_k\,,$$
then after some cancellations, Corollary \ref{C1} reduces to
$$K_k=(3k/4-1)K_{k-1}+\sum_{j=1}^{k-1}K_jK_{k-j}$$
with $K_0=-1/2$. This is the recursive equation found by Tak\'acs, which in turn yields  a Riccati type ODE by considering the exponential generating function of the $K_k$; see Flajolet {\it et al.} \cite{FPV}. We mention the existence of other recursive formulas for the moments of 
$A_{\rm Exc}$, see in particular \cite{louchard2} and the discussion in \cite{janson}. 

Similar calculations apply when more generally the dislocation measure is of beta-type, i.e.
is binary with
$$\nu(\d x) = c x^{\beta}(1-x)^{\beta}\d x\,, \qquad 1/2< x < 1$$
for some $-2<\beta<-1$. These beta-splitting measures have appeared in works of Aldous
\cite{Aldous} on cladograms; see also Section 5.1 in \cite{HMPW}. 
One obtains
\begin{eqnarray*}
\frac{2}{c}a_{k}\, &=&\, \int_0^1 \left(1- x^{k(1-\alpha)}-(1-x)^{k(1-\alpha)}\right) x^{\beta} (1-x)^{\beta}\d x  \\
&=& \, {\rm B}(\beta+1,\beta+1)-2{\rm B}(\beta+1+k(1-\alpha),\beta+1)\\\
&=& \frac{2(2\beta+3)}{\beta+1}{\rm B}(\beta+2,\beta+2)- 2\frac{2\beta+2+k(1-\alpha)}{\beta+1}
{\rm B}(\beta+1+k(1-\alpha), \beta +2)
\end{eqnarray*}
where ${\rm B}(a,b)=\Gamma(a)\Gamma(b)/\Gamma(a+b)$ is the beta function,
and 
\begin{eqnarray*}
\frac{2}{c} a_{j,k}&=& \left(\begin{matrix}
k\\j \end{matrix} \right)\int_0^1 x^{\beta+j(1-\alpha)}(1-x)^{\beta+(k-j)(1-\alpha)}\d x\\
&=& \left(\begin{matrix}
k\\j \end{matrix} \right) {\rm B}(\beta+j(1-\alpha)+1, \beta+(k-j)(1-\alpha)+1)\,.
\end{eqnarray*}

Finally, note that we can also deal with linear combinations of the beta dislocation measures,
which covers for instance the case of Ford's alpha model; see Section 5.2 in \cite{HMPW}. 
\end{section}

\begin{section}{Proof of Theorem \ref{T1}}
This section is devoted to the proof of Theorem \ref{T1}; it relies in four main steps.
In the first sub-section, we establish a priori bounds for the moments of the area, relying on known
properties of the so-called tagged fragment. In the second sub-section, we prove the equation \eqref{E0} in the special case when the dislocation measure is finite. In the third sub-section, we provide the proof of Theorem \ref{T1} by approximation, taking for granted a weak convergence result for the area
that will be established in the final sub-section. 

\subsection{Bounds for the moments of the area}

The purpose of this subsection is to establish
some a priori bounds on the moments $M_k=\E(A^k)$ of the area. Recall the notation \eqref{E8}.

\begin{lemma}\label{L1} We have
$$M_1=1/\Phi(-\alpha)\,$$
and for $k\geq 1$
$$M_k\leq k\,  \frac{k!}{\Phi(-\alpha)\cdots \Phi(-k\alpha)}\,.$$
As a consequence $\E(\exp( c A))<\infty$ whenever $c<\Phi(\infty)$, and in particular the law $\eta$ of $A$ is determined by its entire moments.

\end{lemma} 

\proof It is convenient to work in the setting of interval-fragmentation, i.e. when the fragmentation $\X$ describes the
ranked sequence of the lengths of the interval components of nested open subsets $(\theta(t), t\geq 0)$. Recall from the introduction that $\theta(t)$ can be expressed in the form
$\theta(t)=\{u\in[0,1]: F(u)>t\}$ and note that for every integer $k\geq 1$, there is the identity
$$A^k=\int_0^1 \d u_1 \ldots \int_0^1\d u_k F(u_1)\cdots F(u_k)\,.$$
In other words, we have
$$\E(A^k)=\E(F(U_1)\cdots F(U_k))$$
where $U_1, \ldots, U_k$ are i.i.d. uniform variables on $[0,1]$. This yields
$$M_1=\E(F(U))\quad \hbox{and} \quad M_k\leq k \E(F(U)^k)$$
where $U$ has the uniform distribution on $[0,1]$. 

The variable $F(U)$ should be viewed as the lifetime of the tagged-fragment, i.e. it is the first instant $t$ when the size 
$\chi(t)$ of the interval component of $\theta(t)$ that contains the randomly tagged point $U$ reaches the absorbing state $0$. This variable has the distribution of an exponential functional, 
$$F(U) \stackrel{\mbox{\tiny (law)}}{=}  I = \int_0^{\infty} \exp(\alpha \xi_t)\d t$$
where $\xi=(\xi_t, t\geq 0)$ is a subordinator with Laplace exponent $\Phi$; see Corollary 2 in \cite{Be1}.
Since it is well-known that 
$$\E(I^k)=\frac{k!}{\Phi(-\alpha)\cdots \Phi(-k\alpha)}$$
(cf. for instance Theorem 2 in \cite{BY}), the first two claims are proved, and the last ones follow immediately as the function $\Phi$ increases. \QED

\subsection{The case with finite dislocation rates}

In this subsection, we assume that the fragmentation process has a finite dislocation measure, i.e. $\nu(\m)\in(0,\infty)$.
This means that under $\P$, the process $\X$ stays in state $(1, 0, \ldots)$
during an exponential time $T$ with parameter $\nu(\m)$, and then, independently of the waiting time $T$,  jumps at some
random mass partition $\X(T)$ whose
distribution is given by the normalized dislocation measure $\nu/\nu(\m)$.
We stress that the jump times of $\X$ may nonetheless accumulate right after $T$; in particular then
$\X$ is not a continuous-time Markov chain. Indeed, the first dislocation may produce fragments of arbitrarily small sizes,
which then split again almost instantaneously by self-similarity.

The proof of the following weaker version of Theorem \ref{T1} in this setting is straightforward.

\begin{proposition}\label{P1} Assume that $\X$ has no erosion and finite dislocation measure $\nu$.
Then for every ${\mathcal C}^1$ function $f: \R_+\to \R$ such that the derivative $f'$
has a finite limit at $\infty$, we have
$$\langle \eta, f'\rangle = \int_{\m} \nu(\d \x) \left( \langle \eta, f\rangle - \langle \eta_{\x}, f\rangle \right)\,.$$
\end{proposition}

\proof An application of the strong Markov property at the first dislocation time $T$ and \eqref{E1}
yields the recursive distributional equation (see the survey  \cite{AB} for much more this topic)
$$A = T + \sum_{i=1}^{\infty} X_i(T)^{1-\alpha}A_i$$
where $(A_i)_{i\in\N}$ is a sequence of i.i.d. copies of $A$ which is further independent of 
$\X(T)$. 
This entails
$$\ell(q) = \frac{1}{\nu(\m)+q} \int_{\m}\nu(\d \x) \prod_{i=1}^{\infty}\ell(x_i^{1-\alpha}q)\,,\qquad q\geq 0$$
where $\ell(q)=\E(\exp(-qA))$ is Laplace transform of the area.  By rearrangement, we arrive at
$$-q\ell(q) = \int_{\m}\nu(\d \x) \left( \ell(q)-\prod_{i=1}^{\infty}\ell(x_i^{1-\alpha}q)\right)\,,$$
which is the equation in the statement specified for $f(a)=\e^{-qa}$. This establishes our claim by a standard application of the Stone-Weierstrass theorem (recall that that $A$ has finite moments). \QED

Let us briefly discuss the elementary example when $\nu=\delta_{(1/2,1/2, 0, \ldots)}$. 
We thus start with a single fragment of unit size which splits in two fragments each of size $1/2$
after an exponential time with parameter $1$, and so on. It should be plain that the area can then be expressed in the form
$$A= {\bf e}_{0,1} + 2^{\alpha-1}({\bf e}_{1,1}+{\bf e}_{1,2}) + 2^{2(\alpha -1)}({\bf e}_{2,1}+{\bf e}_{2,2}+{\bf e}_{2,3}+{\bf e}_{2,4}) + \cdots\,,$$
where the ${\bf e}_{i,j}$ for  $j=1, \ldots , 2^i$ and $i=0,1, \ldots$ are i.i.d. standard exponential variables. The Laplace transform of $A$ is thus given by
$$\ell(q) = \prod_{n=0}^{\infty} \left(\frac{1}{1+2^{n(\alpha-1)}q}\right)^{2^n}\,,$$
and the equation 
$$-q\ell(q)=\ell(q)-\ell(2^{\alpha-1}q)^2$$
provided by Proposition \ref{P1} can be checked directly. 

\subsection{Proof of Theorem \ref{T1} by approximation}
In this subsection, we shall derive Theorem \ref{T1}  by approximation from the case when the dislocation measure is finite, taking for granted the weak convergence of the corresponding areas.
Specifically, we introduce the finite measures
$$\nu^{(n)}(\d \x)={\bf 1}_{\{1-x_1>1/n\}}\nu(\d \x)\,,\qquad \x \in\m\,,$$
where $n$ is a sufficiently large integer so that $\nu^{(n)}\not\equiv 0$.
We write $A^{(n)}$ for the area of a self-similar fragmentation process with index $\alpha$, dislocation measure $\nu^{(n)}$ and without erosion, and denote by $\eta^{(n)}$ the distribution of $A^{(n)}$.
Recall \eqref{E1} and set for a generic mass-partition $\x$
$$\eta^{(n)}_{\x}(\d a)= \P\left(\sum_{i=1}^{\infty} x_i^{1-\alpha}A^{(n)}_i\in \d a \right)$$
where $(A^{(n)}_i: i\in\N)$ is a sequence of i.i.d. copies of $A^{(n)}$.
The following crucial lemma will be established in the next sub-section.

\begin{lemma}\label{L2} 
The sequence $(\eta^{(n)}, n\in\N)$ converges weakly to $\eta$ as $n\to\infty$. 
\end{lemma}
The next step to the proof of Theorem \ref{T1} is the following technical result.

\begin{lemma}\label{L3} 
\noindent {\rm (i)} Let  $f:\R_+\to \R$ be continuous and bounded. Then for every $\x\in\m$, 
$$ \lim_{n\to\infty} 
\langle \eta^{(n)}_{\x},f\rangle =   \langle \eta_{\x},f\rangle\,.$$

\noindent {\rm (ii)} Let  $f:\R_+\to \R$ be a ${\mathcal C}^1$ function with
$f'(y)=O(y^p)$ as $y\to\infty$  for some $p>0$, and set $\|f'\|=\sup \{|f'(x)|/(1+x)^p: x\geq 0\}$.
There is a constant $c$ depending only on $p$ and the characteristics of the fragmentation 
such that for every mass-partition $\x$ and every $n$ 
$$|\langle \eta^{(n)},f\rangle- \langle \eta^{(n)}_{\x},f\rangle |\leq c \| f'\|(1-x_1)\,.$$
\end{lemma}

\proof (i) 
Denote by $\kappa$ the cumulant of $A$, i.e.
$\E(\exp(-qA))=\exp(-\kappa(q))$ for $q\geq 0$, and by $\kappa^{(n)}$ that of $A^{(n)}$. 
If $(A^{(n)}_i, i\in\N)$ is a sequence of i.i.d. copies of $A^{(n)}$,  then we have
$$\E\left(\exp\left(-q \sum_{i=1}^{\infty} x_i^{1-\alpha} A^{(n)}_i\right)\right)=\exp\left(-\sum_{i=1}^{\infty}
\kappa^{(n)}(x_i^{1-\alpha}q)\right)\,.$$
We know from Lemma \ref{L2}  that $\lim_{n\to\infty}\kappa^{(n)}(x_i^{1-\alpha}q)=\kappa(x_i^{1-\alpha} q)$ for every $q\geq 0$ and $i\in\N$.
Further, $\kappa^{(n)}$ is a concave increasing function with $\kappa^{(n)}(0)=0$, and its derivative at $0$
is given by $\E(A^{(n)})$. Recall also from Lemma \ref{L1} that
$\E(A^{(n)})=1/\Phi^{(n)}(-\alpha)$ where
$$\Phi^{(n)}(q)=\int_{\m}\left(1-\sum_{i=1}^{\infty}x_i^{1+q}\right)\nu^{(n)}(\d \x)
=\int_{\m}{\bf 1}_{\{1-x_1>1/n\}}\left(1-\sum_{i=1}^{\infty}x_i^{1+q}\right)\nu(\d \x)\,.$$
Plainly the sequence $\Phi^{(n)}(-\alpha)$ increases 
and there are thus the bounds
$$\kappa^{(n)}(x_i^{1-\alpha}q)\leq x_i^{1-\alpha}q/\Phi^{(n)}(-\alpha) \leq  c x_i^{1-\alpha}q \,.$$
Since $\sum_{i=1}^{\infty}x_i^{1-\alpha} \leq 1$, we conclude by dominated convergence that
$$\lim_{n\to\infty} \sum_{i=1}^{\infty}
\kappa^{(n)}(x_i^{1-\alpha}q)= \sum_{i=1}^{\infty}
\kappa(x_i^{1-\alpha}q)\,.$$
By Laplace inversion, this yields our claim.

(ii) 
Recall that  $A^{(n)}_1, \ldots$ are  i.i.d. copies of $A^{(n)}$,  and set
$S^{(n)}_{\x}=\sum_{i=2}^{\infty} x_i^{1-\alpha} A^{(n)}_i$ where $\x=(x_1, \ldots)$ is a generic mass-partition.
We have
\begin{eqnarray*}
| \langle \eta^{(n)}, f\rangle - \langle \eta^{(n)}_{\x}, f\rangle| 
&\leq& \E\left(\left| f(A_1^{(n)})-f\left(x_1^{1-\alpha} A_1^{(n)} + S^{(n)}_{\x}\right)\right|\right)\\
&\leq & \|f\| \E\left((1+A_1^{(n)}+S^{(n)}_{\x})^p\left((1-x_1^{1-\alpha})A_1^{(n)} + S^{(n)}_{\x}\right) \right)\,.
\end{eqnarray*}

Recall that $\sum_1^{\infty}x_i\leq 1$. 
Applying Rosenthal's inequality  for the moments of sum of independent nonnegative variables and then Jensen's inequality, we obtain
$$\E\left((1+A_1^{(n)}+S^{(n)}_{\x})^{2p}\right)\leq c(2p)
(1+M_{2p}^{(n)})$$  
and 
$$\E\left(\left((1-x_1^{1-\alpha})A_1^{(n)} + S^{(n)}_{\x}\right) ^2\right)\leq c(2)(1-x_1)^2M_2^{(n)}\,,$$
where $c(2p)$ and $c(2)$ are some numerical constants and $M^{(n)}_k$ denotes the $k$-th moment of $A^{(n)}$. We conclude from
 H\"older's inequality that
$$| \langle \eta^{(n)}, f\rangle - \eta^{(n)}_{\x}, f\rangle|  \leq c\| f \| (1-x_1)\sqrt{M\n_2(1+ M_{2p}^{(n)})}\,.$$
We can complete the proof with an appeal to Lemma \ref{L1},
recalling that the sequence
$n\mapsto \Phi^{(n)}(q)$ 
increases for every $q>0$ and thus $\sup_n\sqrt{M\n_2(1+ M_{2p}^{(n)})}<\infty$. 
 \QED

We can now proceed to the proof of Theorem \ref{T1}.

\noindent {\bf Proof of Theorem \ref{T1}: \hskip8pt}  
We first suppose that  $f$ and $f'$ are bounded.
Then \eqref{E0} follows from Proposition \ref{P1} and Lemma \ref{L3} by dominated convergence. 
Next, we only assume that $f$ is of class ${\mathcal C}^1$ with $f'(y)=O(y^p)$ as $y\to\infty$. 
Then it is easy to construct a sequence $(f_n: n\in\N)$ of functions  of class ${\mathcal C}^1$ with  $f_n$ and $f'_n$ bounded such that $f_n\to f$ and $f'_n\to f'$ pointwise  and $\sup_n \| f'_n\|<\infty$ in the sense of Lemma \ref{L2}. This implies that we also have $\sup_n\sup_{x\geq 0}|f_n(x)/(1+x^{p+1})|)<\infty$.

Since all moments of $A$ are finite,
we deduce from Lemmas \ref{L2} and \ref{L3}(i) that 
$$\langle \eta, f_n\rangle \to \langle \eta, f\rangle\ ,\ 
\langle \eta, f'_n\rangle \to \langle \eta, f'\rangle\ \hbox{ and } \ 
\langle \eta_{\x}, f_n\rangle \to \langle \eta_{\x}, f\rangle$$
where $\x$ is an arbitrary mass-partition. 
Since $\sup_n \| f'_n\|<\infty$, Lemma \ref{L3}(ii) enables us to apply dominated convergence,
and we conclude that \eqref{E0} holds.

We now turn our attention to uniqueness; we consider an arbitrary solution $\eta'$
to \eqref{E0} and write $M'_k$ for the $k$-th moment of $\eta'$.
We have already observed in the introduction that  the first moment $M'_1$  of $\eta'$  can be computed in terms of the function $\Phi$ defined by \eqref{E8}. More generally, specifying \eqref{E0} for $f(x)=x^k$
yields an equation of the form
$$\Phi(k(1-\alpha)-1)M'_k= \Psi_k(M'_1,\ldots, M'_{k-1})$$
for a certain multinomial function $\Psi_k$. Hence $M'_k=M_k$ for every $k\in\N$,
and since we know from Lemma \ref{L1} that the moment problem for $\eta$ is determined,
this concludes the proof. \QED

\subsection{Proof of Lemma \ref{L2}}
We shall finally establish Lemma \ref{L2} using the framework of homogeneous fragmentations
with values  in the space $\p_{\N}$ of partitions of $\N$. 
Given a measure $\nu$ which fulfills \eqref{E9}, we first consider a
homogeneous fragmentation $\Pi=(\Pi(t), t\geq 0)$ 
with dislocation measure $\nu$ and no erosion.
For every $t\geq 0$ the random partition $\Pi(t)$ is exchangeable and we write $|\Pi_i(t)|$ for the asymptotic frequency of the $i$-th block of  $\Pi(t)$. 
The self-similar fragmentation $\X$ is related to  $\Pi$ by a sort of time-change described in Theorem  
3.3 of \cite{RFCP}, and if we  introduce 
$$\a = \int_0^{\infty}\left( \sum_{i=1}^{\infty}|\Pi_i(t)|^{1-\alpha}\right) \d t\,,$$
then the connexion between $\Pi$ and $\X$ implies that $\a$ and $A$ have the same distribution.

It will be convenient to approximate $\a$ by Riemann sums. More precisely, for every integer $k\geq 1$ we define
$$\a_k=\frac{1}{k}\sum_{\ell=1}^{k^2} \sum_{i=1}^{\infty}|\Pi_i(\ell/k)|^{1-\alpha}\,.$$
The $L^1$ distance between these two quantities is easily computed in terms of the function $\Phi$ defined in \eqref{E8}. 

\begin{lemma}\label{L4} For every $k\geq 1$, we have
$$\E( | \a-\a_k|)= \frac{1}{\Phi(-\alpha)} - k^{-1}\frac{1-\exp(-k\Phi(-\alpha))}
{\exp(\Phi(-\alpha)/k)-1}\,.$$
\end{lemma}

\proof  The process $|\Pi(\cdot)|^{\downarrow}$ of the ranked sequence of the asymptotic frequencies of $\Pi(\cdot)$  is a mass-fragmentation; as a consequence $t\mapsto \sum_{i=1}^{\infty}|\Pi_i(t)|^{1-\alpha}$ decreases and  we have $\a_k\leq \a$. Hence $\E( | \a-\a_k|)=\E(\a)-\E(\a_k)$
and the stated formula now follows from the fact that
$$\E\left( \sum_{i=1}^{\infty}|\Pi_i(t)|^{1-\alpha}\right)= \exp(-t\Phi(-\alpha))\,;$$
see Corollary 2.4(i) and Theorem 3.2 in \cite{RFCP}. \QED

Next, for every integer $n$, we write $\Pi\n$ for the homogeneous  fragmentation with dislocation measure $\nu\n(\d \x)={\bf 1}_{\{1-x_1>1/n\}}\nu(\d \x)$ and no erosion.
 We point at the following weak convergence.

\begin{lemma}\label{L5} The sequence of processes $(\Pi\n: n\in\N)$ converges
in the sense of finite dimensional distributions to $\Pi$ as $n\to\infty$.
\end{lemma}

\proof  If $\varphi\in\p_{\N}$  is a partition and $k\geq 1$ an integer, then we denote by
 $\varphi_{\mid [k]}$ the restriction of $\varphi$ to the set of the $k$ first integers,  $[k]=\{1, \ldots, k\}$.
 We also endow  $\p_{\N}$  with the ultra-metric 
 $$d(\varphi, \varphi')=1/\sup\{k\geq 1: \varphi_{\mid [k]} = \varphi'_{\mid [k]}\}\,;$$
  cf. Lemma 2.6 in \cite{RFCP}.

The restriction $\Pi_{\mid [k]}(t)$  of $\Pi(t)$ to $[k]$ is a Markov chain in continuous time, and we have to verify that $\Pi\n_{\mid [k]}$ converges
in the sense of finite dimensional distributions to $\Pi_{\mid [k]}$ as $n\to\infty$, for each $k$.
This is equivalent to checking the convergence of the corresponding jump rates of the Markov chains.

For every non-trivial partition $\gamma$ of $[k]$,  we write
 $$q_{\gamma}=\lim_{t\to 0}t^{-1}\P(\Pi_{\mid [k]}(t)=\gamma)\,.$$
Recall from Theorem 3.1 and Proposition 3.2 in \cite{RFCP} that the jump rate $q_{\gamma}$ can be expressed in terms of the dislocation measure $\nu$ as 
 $$q_{\gamma} = \int_{\m} \varrho_{\x}(\gamma) \nu(\d \x)$$
 where $\varrho_{\x}(\gamma)$ is the distribution of Kingman's paintbox process $\pi_{\x}$ based on $\x$. 
 This means that we consider a sequence $\xi_1, \ldots, \xi_k $ of i.i.d. variables with
 $\P(\xi_1=i)=x_i$ for $i\geq 1$ and $\P(\xi_1=0)=1-\sum_1^{\infty}x_i$,
 and  $\pi_{\x}$ is the exchangeable random partition which is obtained by declaring that two integers  $i\neq j$ are in the same block of $\pi_{\x}$ if and only if $\xi_i= \xi_j\neq 0$. 
 
 Writing   $q_{\gamma}\n$ for the jump rate of $\Pi\n$, we thus have
 $$  q_{\gamma}-q_{\gamma}\n=\int_{\m}{\bf 1}_{\{1-x_1\leq 1/n\}}
 \varrho_{\x}(\gamma) \nu(\d\x)\,.$$
It is plain from the paintbox construction that for every mass-partition $\x$,
the probability that the paintbox process based on $\x$ yields the trivial partition on $[k]$
is at least $x_1^k$. Thus $\varrho_{\x}(\gamma)\leq 1-x_1^k \leq k(1-x_1)$
for every non-trivial partition $\gamma$ of $[k]$.
We conclude from \eqref{E9} that $\lim_{n\to\infty}q_{\gamma}\n=q_{\gamma}$, which establishes our claim. \QED

We are now in shape to prove Lemma \ref{L2}

\noindent {\bf Proof of Lemma \ref{L2}:}Ê\hskip8pt
The space $\m$ of mass-partitions is a compact metric space when endowed with the uniform distance (Proposition 2.1 in \cite{RFCP}) and the map $\x\to \sum_{i=1}^{\infty} x_i^{1-\alpha}$
is continuous and bounded. 
We write $|\Pi(t)|^{\downarrow}$ for the sequence of the asymptotic frequency of  $\Pi(t)$ ranked in the decreasing order.
Recall from Proposition 2.9 in \cite{RFCP} that Lemma \ref{L5} entails the convergence in the sense of finite dimensional distributions of $(|\Pi\n(t)|^{\downarrow}, t\geq 0)$ towards  $(|\Pi(t)|^{\downarrow}, t\geq 0)$. Therefore, for every 
$k\geq 1$, we have also in the obvious notation
\begin{equation} \label{E7}
\lim_{n\to\infty} \a\n_k= \a_k \quad\hbox{in law.}
\end{equation}
 
Let $f:\R_+\to\R$ be a bounded function which is globally Lipschitz-continuous, so by the triangle inequality
$$|\E(f(\a))-\E(f(\a\n))|\leq |\E(f(\a_k))-\E(f(\a\n_k))| + c_f (\E(|\a-\a_k|)+  \E(|\a\n-\a\n_k|))\,.$$
It follows readily from Lemma \ref{L4} that
$$\lim_{k\to\infty} \E(|\a\n-\a\n_k|)=0\qquad \hbox{uniformly in } n\,,$$
so for every $\varepsilon>0$, we can find an integer $k$ sufficiently large
such that 
$$\E(|\a-\a_k|)+  \E(|\a\n-\a\n_k|)\leq \varepsilon/2c_f \qquad \hbox{ for all }n\,,$$ 
and then we use 
\eqref{E7} to find an integer $n_{\varepsilon}$ such that $ |\E(f(\a_k))-\E(f(\a\n_k))| \leq \varepsilon/2$
whenever $n\geq n_{\varepsilon}$.
\QED
 \end{section}

\vskip 5mm 
\noindent{\bf Acknowledgments}.  This work has been supported by ANR-08-BLAN-0220-01.

  \end{document}